\documentclass[twoside]{amsart}
\usepackage{latexsym}
\usepackage{amssymb,amsmath,amsopn}
\usepackage[dvips]{graphicx}   
\usepackage{color,epsfig}      
\usepackage{tikz} 
\usepackage{subcaption} 
\usepackage{eurosym}

\renewcommand{\Re}{\mathbb R}
\newcommand{\Z}{\mathbb Z}
\newcommand{\C}{\mathbb C}
\newcommand{\B}{\mathbf B}

\newcommand{\Sph}{\mathbb{S}}

\DeclareMathOperator{\bd}{bd}

\DeclareMathOperator{\conv}{conv}

\DeclareMathOperator{\tr}{tr}

\DeclareMathOperator{\aff}{aff}

\DeclareMathOperator{\Aut}{Aut}

\DeclareMathOperator{\Id}{Id}
\DeclareMathOperator{\lcm}{lcm}

\newtheorem{theorem}{Theorem}
\newtheorem{prop}{Proposition}[section]

\newtheorem{lemma}[prop]{Lemma}
\newtheorem{problem}[prop]{Problem}

\theoremstyle{definition}
\newtheorem{remark}[prop]{Remark}

\begin{document}

\title[Affinely regular polygons]{A characterization of affinely regular polygons}
\author[Z. L\'angi]{Zsolt L\'angi}
\address{MTA-BME Morphodynamics Research Group and Dept.\ of Geometry, Budapest University of
Technology and Economics\\
Egry J\'ozsef u. 1.\\
1111 Budapest\\
Hungary, 1111}
\email{zlangi@math.bme.hu}

\begin{abstract}
In 1970, Coxeter gave a short and elegant geometric proof showing that if $p_1, p_2, \ldots, p_n$ are vertices of an $n$-gon $P$ in cyclic order, then $P$ is affinely regular if, and only if there is some $\lambda \geq 0$ such that $p_{j+2}-p_{j-1} = \lambda (p_{j+1}-p_j)$ for $j=1,2,\ldots, n$. The aim of this paper is to examine the properties of polygons whose vertices $p_1,p_2,\ldots,p_n \in \C$ satisfy
the property that $p_{j+m_1}-p_{j+m_2} = w (p_{j+k}-p_j)$ for some $w \in \C$ and $m_1,m_2,k \in \Z$.
In particular, we show that in `most' cases this implies that the polygon is affinely regular, but in some special cases there are polygons which satisfy this property but are not affinely regular. The proofs are based on the use of linear algebraic and number theoretic tools.
In addition, we apply our method to characterize polytopes with certain symmetry groups.
\end{abstract}

\thanks{Partially supported by the National Research, Development and Innovation Office, NKFI-119670}

\maketitle


\section{Introduction}

In recent times affinely regular polygons have played an important role in various mathematical problems. They appear, for example, in geometry as extremal sets in optimization problems \cite{GronchiLonginetti} or in the famous Napoleon-Barlotti Theorem, redisovered by Gerber in 1980 \cite{Gerber}, or in discrete tomography \cite{Gardner}; in linear algebra as sets of vectors cyclically permuted by unimodular matrices.
In the 1970s, Coxeter started a systematic investigation of the properties of affinely regular polygons \cite{Coxeter1} in the Euclidean plane, which was later continued both in the Euclidean plane \cite{Coxeter2}, \cite{Huck}, \cite{Nicollier}, and in affine planes in general \cite{CraatsSimonis}, \cite{Kiss}, \cite{KorchmarosSzonyi}. For a survey about these properties, particularly about characterizations of affinely regular polygons, the interested reader is referred to \cite{FischerJamison}.

The following theorem appeared in \cite{FischerJamison}, where the authors attribute it to Coxeter \cite{Coxeter1}.

\begin{theorem}[Coxeter]\label{thm:coxeter}
Let $p_1, p_2, \ldots, p_n$ be the vertices of an $n$-gon $P$ in the Euclidean plane $\Re^2$, in cyclic order. If there is some real number $\lambda \geq 0$ such that $p_{j+2}-p_{j-1} = \lambda (p_{j+1}-p_j)$ for $j=1,2,\ldots,n$, then $P$ is affinely regular.
\end{theorem}

Our aim is to examine the following, more general problem, where, for simplicity, we identify the Euclidean plane $\Re^2$ by the complex plane $\C$.

\begin{problem}\label{prob:main}
Find conditions on the integers $m_1,m_2,k \in \Z$, and the complex number $w \in \C$ such that if the vertices $p_1,p_2,\ldots,p_n \in \C$ of a polygon $P$ satisfy the equations $p_{j+m_1}-p_{j+m_2} = w (p_{j+k}-p_j)$ for $j=1,2,\ldots,n$, then $P$ is affinely regular.
\end{problem}

Before stating our main result, we introduce some definitions and notation.
For $n \geq 4$, by an \emph{$n$-gon} $P = (p_1,p_2,\ldots,p_n)$, we mean an (ordered) $n$-tuple of points in $\C$, where the indices of the points are understood modulo $n$. We call the points $p_1,p_2,\ldots,p_n$ the \emph{vertices} of $P$. To avoid degenerate cases, throughout the paper we only deal with $n$-gons whose vertices are pairwise distinct. We say that $P=(p_1,p_2,\ldots,p_n)$ is \emph{affinely regular} (respectively, \emph{regular}), if there is an affine transformation (respectively, isometry) $\phi$ such that $\phi(p_j)=p_{j+1}$ for $j=1,2,\ldots,n$. In particular, if we set $\epsilon = \cos \frac{2\pi}{n}+i \sin \frac{2\pi}{n}$, then the polygon $(1,\epsilon^t,\epsilon^{2t}, \ldots, \epsilon^{((n-1)t})$ is a regular $n$-gon for every value of $t$ relatively prime to $n$.

For any nonzero integers $k_1,k_2,\ldots,k_s \in \Z$, we denote by the symbol $\gcd(k_1,k_2,\ldots,k_s)$ the greatest common divisor of $k_1,k_2,\ldots,k_s$, i.e. the largest po\-si\-ti\-ve integer $x$ satisfying $x \mid k_j$ for all $j=1,2,\ldots, s$. For the least common multiple of the nonzero integers $k_1,k_2,\ldots,k_s \in \Z$, we use the notation $\lcm (k_1,k_2,\ldots,k_s)$.

Our main result is as follows.

\begin{theorem}\label{thm:main}
Let $n,m_1,m_2,k \in \Z$ be integers such that $n \geq 4$,  neither $k$ nor $m_1-m_2$ is divisible by $n$, and $\gcd(n,k,m_1-m_2) = 1$. If $n$ is even, assume that $n > 2 \gcd(n,k) \gcd(n,m_1-m_2)$. Let $P = (p_1,p_2,\ldots,p_n)$ be an $n$-gon whose vertices satisfy the equations $p_{j+m_1}-p_{j+m_2} = w (p_{j+k}-p_j)$ for $j=1,2,\ldots,n$ and for some $w \in \C$ with $|w| \neq 1$ and independent of $j$.
\begin{enumerate}
\item[(\ref{thm:main}.1)] If $w \notin \Re$, then $P$ is a regular $n$-gon, and $m_1+m_2 \not\equiv k \pmod n$.
\item[(\ref{thm:main}.2)] If $w \in \Re$, then $P$ is an affinely regular $n$-gon, and $m_1 + m_2 \equiv k \pmod n$.
\end{enumerate}
\end{theorem}

As an application of our method, we also investigate another problem.
Schulte and Williams~\cite{SchulteWilliams} proved that every finite group is the combinatorial automorphism group of a suitably chosen convex polytope. This result was strengthened by Doignon \cite{Doignon} who proved the existence of such a convex polytope with the additional property that the combinatorial and the geometric automorphism groups of the polytope are equal.
Thus, it is a natural question to ask the following: given a finite group $\Gamma$, which convex polytopes have geometric automorphism groups containing $\Gamma$ as a subgroup. The special case where we look for $d$-dimensional polytopes with $d+3$ vertices and $\Gamma = D_{d+3}$ appeared in \cite{GHorvathLangi} as Problem 3. We solve a more general problem: we characterize the $d$-dimensional polytopes with $n$ vertices whose geometric automorphism group contains $D_n$ as a subgroup. This result can be regarded as a first step towards solving Problem 1 \cite{GHorvathLangi}, and we note that $D_n$ is the combinatorial automorphism group of a cyclic polytope in an even dimensional space.

Before stating our second result, we introduce some additional notation.
Let $d \geq 2$, $n > d$, $s=\lfloor d/2 \rfloor$, $0 < k_1 < k_2 < \ldots < k_s < \frac{n}{2}$, and for $m=1,2,\ldots,n$, let
\[
q_m = \frac{1}{\sqrt{s}} \left( \cos \frac{2 k_1 m \pi}{n}, \sin \frac{2k_1 m \pi}{n}, \ldots, \cos \frac{m k_s\pi}{n}, \sin \frac{2m k_s\pi}{n} \right) \in \Re^d
\]
if $d$ is even, and
\[
q_m = \frac{1}{\sqrt{s+1}} \left( \cos \frac{2 k_1 m \pi}{n}, \sin \frac{2k_1 m \pi}{n}, \ldots,
\sin \frac{2m k_s\pi}{n}, (-1)^m \right) \in \Re^d
\]
if $d$ is odd. Furthermore, set $Q(k_1,k_2,\ldots,k_s) = \conv \{ q_1,q_2,\ldots,q_n\}$. Note that as $|q_m| = 1$ for every value of $m$, the vertices of $Q(k_1,k_2,\ldots,k_s)$ are $q_1,q_2,\ldots,q_n$.
Finally, for a convex polytope $P$ in $\Re^d$, let $\Aut(P)$ denote the geometric automorphism group of $P$; that is, the group of isometries that leave $P$ invariant.
Our second result is the following.

\begin{theorem}\label{thm:symmetric}
Let $d \geq 2$, and $P \subset \Re^d$ be a $d$-dimensional convex polytope with vertices $p_1,p_2,\ldots, p_n$, where $n \geq 5$, and $n > d$. Then the following are equivalent.
\begin{enumerate}
\item{(\ref{thm:symmetric}.1)} For $j,k=1,2,\ldots,n$, the value of $|p_{j+k}-p_j|$ is independent of the value of $j$.
\item{(\ref{thm:symmetric}.2)} There is some $\phi \in \Aut(P)$ such that $\phi(p_j) = p_{j+1}$ for $j = 1,2,\ldots,n$.
\item{(\ref{thm:symmetric}.3)} There are some $0 < k_1 < k_2 < \ldots < k_{\lfloor d/2 \rfloor} < \frac{n}{2}$ such that
$P$ is similar to $Q(k_1,k_2,\ldots,k_{\lfloor d/2 \rfloor})$.
\end{enumerate}
\end{theorem}

In Section~\ref{sec:affinelyregular} we prove Theorem~\ref{thm:main}. The proof is based on Lemma~\ref{lem:conway} about the ratios of the lengths of the diagonals and sides of a regular polygon, which we prove in Section~\ref{sec:conway}. This lemma might be of independent interest. In Section~\ref{sec:polytopes} we present the proof of Theorem~\ref{thm:symmetric}. Finally, in Section~\ref{sec:remarks} we collect our remarks. In this section, in particular, we give examples showing why the conditions in Theorem~\ref{thm:main} are necessary.

\section{Proof of Theorem~\ref{thm:main}}\label{sec:affinelyregular}

During the proof we let $m = m_2-m_1$, and $\epsilon = \cos \frac{2\pi}{n} + i \sin \frac{2\pi}{n}$. Furthermore, since the vertices of $P$ are pairwise distinct and neither $k$ nor $m$ is divisible by $n$, we observe that $w \neq 0$. For simplicity, we assume that $0 < k,m_1,m_2 < n$.

Since the points $p_j$, $j=0,1,\ldots,n-1$ satisfy the conditions $p_{j+m_1}-p_{j+m_2} = w (p_{j+k}-p_j)$ for every value of $j$, the solutions for $P$ are exactly the complex solutions of the homogeneous system of the linear equations $w p_j - w p_{j+k}+ p_{j+m_1}-p_{j+m_2} = 0$, where $j=0,1,\ldots, n-1$.
In other words, they are the eigenvectors in $\C^n$ of the coefficient matrix $C$ of this system corresponding to the eigenvalue $0$.
The coefficient matrix $C$ of this system is a circulant matrix; i.e., a square matrix with the property that its $(j+1)$st row is the shift of its $j$th row by one to the right.

The eigenvalues and eigenvectors of a circulant matrix are known \cite{Davis}: If the first row of $C$ is the vector $c=(c_0,c_1,\ldots,c_{n-1})$, then its eigenvalues are  $\mu_t = \sum_{j=0}^{n-1} c_j \epsilon^{jt}$, where $t=0,1,\ldots,n-1$, with the corresponding eigenvector $v_t = \left( \epsilon^t, \epsilon^{2t}, \ldots,\epsilon^{(n-1)t}\right)$.

In the remaining part of the proof we examine for which values of $t$ can $\mu_t$ be equal to zero, depending on the values of $m_1,m_1,k$ and $w$.
During this examination, without loss of generality, we assume that $c_{m_2}=-1$, $c_0 = w$, $c_k = -w$, $c_{m_1}=1$.

First, note that $\mu_0 = -1+w-w+1 = 0$, and thus, $v_0$ is a solution for every value of $w$.

Next we show that if $n$ is even, then $\mu_{n/2} \neq 0$. Indeed, suppose for contradiction that $\mu_{n/2} = 0$.
Then we have
\[
w \left( 1-(-1)^k\right) = (-1)^{m_1} \left( (-1)^{m} - 1 \right).
\]
Since $w \neq 0$ and $w \neq \pm 1$, this equality implies that $1-(-1)^k = (-1)^{m} - 1 = 0$; that is, that both $k$ and $m$ are even. In this case $\gcd(n,k,m) \geq 2$; a contradiction, and hence $\mu_{n/2} \neq 0$.

Now we have that there is some $ t \neq 0, \frac{n}{2}$ with $\mu_t = 0$.
Consider such a value $t$. Then $w \left( \epsilon^{tk}-1 \right) = \epsilon^{tm_1}- \epsilon^{tm_2}$. We show that $\epsilon^{tk}-1 \neq 0$.
Indeed, if $\epsilon^{tk}-1 \neq 0$, then $\epsilon^{tm_1}- \epsilon^{tm_2} = 0$, or in other words, $tk \equiv 0$ and $tm_1 \equiv tm_2 \pmod n$. This implies that
$t \equiv 0 \pmod{ \frac{n}{\gcd(k,n)}}$ and $t \equiv 0 \pmod{\frac{n}{\gcd(m,n)}}$, which yields that $t \equiv 0 \pmod{\lcm ( \frac{n}{\gcd(k,n)}, \frac{n}{\gcd(m,n)} ) = n}$; a contradiction.
Thus, we have $w = \frac{\epsilon^{tm_1}- \epsilon^{tm_2}}{\epsilon^{tk}-1}$. Since the complex numbers $1 , \epsilon, \ldots, \epsilon^{n-1}$ are the vertices of a regular $n$-gon, in this case $w$ is the ratio of the lengths of two diagonals or sides of a regular $n$-gon.

Now we examine the case that some $t, t' \neq 0, \frac{n}{2}$ satisfy $\mu_t = \mu_{t'} = 0$. This implies that
\begin{equation}\label{eq:indices}
\frac{\epsilon^{tm_1}- \epsilon^{tm_2}}{\epsilon^{tk}-1} = \frac{\epsilon^{t'm_1}- \epsilon^{t'm_2}}{\epsilon^{t'k}-1} .
\end{equation}

Before proceeding further, observe that since both sides in (\ref{eq:indices}) are equal to $w$, $w \in \Re$ yields that if $\mu_t = 0$ then $\mu_{n-t} = 0$, and if for some $t \neq 0, \frac{n}{2}$ we have $\mu_t = \mu_{n-t} = 0$, then $w \in \Re$.

First, we only investigate the absolute values of the sides in (\ref{eq:indices}), and apply the next lemma, which we prove in Section~\ref{sec:conway}.

\begin{lemma}\label{lem:conway}
Let $n \geq 4$, and let $Q$ be a regular $n$-gon in the Euclidean plane $\Re^2$, with vertices $q_0,q_1,\ldots, q_{n-1}$ in counterclockwise order. For any $1 \leq j \leq \frac{n}{2}$, set $d_j = |q_j-q_0|$. If for some values $1 \neq k,l,k',l' \leq \frac{n}{2}$, we have $1 \neq \frac{d_k}{d_l}=\frac{d_{k'}}{d_{l'}}$, then $k=k'$ and $l=l'$.
\end{lemma}

By Lemma~\ref{lem:conway}, it follows that $t'k \equiv tk \pmod n$ or $t'k \equiv -tk \pmod n$.

\emph{Case 1}, $t'k \equiv tk \pmod n$.
Then $\epsilon^{tk}-1 = \epsilon^{t'k}-1$, which implies that $\epsilon^{tm_1}- \epsilon^{tm_2} = \epsilon^{t'm_1}- \epsilon^{t'm_2}$.
Since both sides are diagonals or sides of a regular $n$-gon, at least one of the following two systems of congruences holds.
\begin{enumerate}
\item[(i)] $t'm_2 \equiv t m_2$ and $t' m_1 \equiv tm_1 \pmod n$.
\item[(ii)] $n$ is even, and $tm_1+\frac{n}{2} \equiv t'm_2$ and $tm_2+\frac{n}{2} \equiv t'm_1 \pmod n$.
\end{enumerate}

In the first case, subtracting the two congruences we obtain $t'm \equiv tm \pmod n$, which, as $\gcd(n,k)$ and $\gcd(n,m)$ are relatively prime, yields that $t' \equiv t \pmod n$.
In the second case, it follows similarly that $t' m \equiv -tm \pmod n$. Thus, $t'$ is a solution of the simultaneous congruence system
\begin{equation}\label{eq:congruences}
t' \equiv t \pmod{\frac{n}{\gcd(n,k)}}, \quad \mathrm{and} \quad t' \equiv -t \pmod{\frac{n}{\gcd(n,m)}}.
\end{equation}
It is well-known \cite{IrelandRosen} that this system has a solution for $t'$ if, and only if $t \equiv -t \pmod{\gcd ( \frac{n}{\gcd(n,k)}, \frac{n}{\gcd(n,m)} ) = \frac{n}{\gcd(n,k)\gcd(n,m)}}$, and in this case its solution is unique $\mod \lcm \left( \frac{n}{\gcd(n,k)}, \frac{n}{\gcd(n,m)} \right) = n$. In this special case the unique solution, if it exists, can also be computed. Indeed, since $\gcd(n,k)$ and $\gcd(n,m)$ are relatively prime, there are integers $x,y$ such that $1 = x \gcd(n,k)+y\gcd(n,m)$. Thus, the congruences in (\ref{eq:congruences}) can be written in the form $t' \gcd(n,k) \equiv t \gcd(n,k)$ and $t' \gcd(n,m) \equiv -t\gcd (n,m) \pmod n$, which implies that
\[
t' \equiv \left( x\gcd (n,k) - y \gcd (n,m) \right) t \pmod n.
\]
Reversing the roles of $t$ and $t'$ we obtain that $t \equiv \left( x \gcd(n,k) - y \gcd(n,m) \right) t' \pmod n$. This implies, in particular, that $\gcd(t,n) = \gcd(t',n)$.

\emph{Case 2}, $t'k \equiv -tk \pmod n$.
Then $\epsilon^{t'k}-1= \epsilon^{-tk}-1 =\frac{-1}{\epsilon^{tk}} \left( \epsilon^{tk}-1\right)$, from which we have that $\epsilon^{t'm_1}-\epsilon^{t'm_2} = \epsilon^{t(m_2+k)-\epsilon^{t(m_1+k)}}$.
Thus, at least one of the following holds.
\begin{enumerate}
\item[(i)] $t'm_2 \equiv t m_1+ tk$ and $t' m_1 \equiv tm_1 +tk\pmod n$.
\item[(ii)] $n$ is even, and $tm_1+tk + \frac{n}{2} \equiv t'm_1$ and $tm_2+tk + \frac{n}{2} \equiv t'm_2 \pmod n$.
\end{enumerate}
In the first case, similarly like in Case 1, we obtain that $t' \equiv -t \pmod n$. In the second case, by the argument used in Case 1, we obtain that
we have a solution if, and only if $t \equiv -t \pmod{\frac{n}{\gcd(n,k)\gcd(n,m)}}$, and in this case the solution is unique $\mod n$. Furthermore, if there is a solution, then $t' \equiv \left( y \gcd(n,m) - x \gcd(n,k) \right) t \pmod n$ and $t \equiv \left( y (n,m) - x(n,k) \right) t' \pmod n$, yielding that $\gcd(t,n) = \gcd (t',n)$.

Summing up, we are left with one of the following.

\begin{itemize}
\item[(a)] We have $\mu_t = 0$ if, and only if $t=0$.
\item[(b)] We have $w \notin \Re$, and there  is exactly one value of $t \neq 0$ such that $\mu_t = 0$.
\item[(c)] We have $w \in \Re$, and there is a unique $t \neq 0, \frac{n}{2}$ such that $\mu_s = 0$ if, and  only if $s \in \{ 0,t,n-t\}$.
\item[(d)] We have $w \notin \Re$ and $2 | n$, and there are some distinct values $t,t' \neq 0, \frac{n}{2}$ such that $t' \neq n-t$, and $\mu_s = 0$ if, and only if $s \in \{ 0,t,t'\}$. In particular, in this case $2t \equiv 0 \pmod{\frac{n}{\gcd(n,k)\gcd(n,m)}}$, and $\gcd(t,n) = \gcd(t',n)$.
\item[(e)] We have $w \in \Re$ and $2 | n$, and there are some distinct values $t,t' \neq 0, \frac{n}{2}$ such that $t' \neq n-t$, and
$\mu_s = 0$ if, and only if $s \in \{ 0,t,t',n-t,n-t'\}$. In particular, in this case $2t \equiv 0 \pmod{\frac{n}{\gcd(n,k)\gcd(n,m)}}$, and $\gcd(t,n) = \gcd(t',n)$.
\end{itemize}

If the conditions in (a) hold, then the only solutions are of the form $P= z v_0 = (z,z,\ldots, z)$ for some point $z \in \C$. Thus, in this case the vertices of $P$ are not pairwise distinct; a contradiction.
We show that under the conditions in (d) or (e), there is no solution as well.
We prove this statement for (d), as for (e) a similar argument can be applied. Assume that there are some distinct values $t,t' \neq 0, \frac{n}{2}$ such that $t' \neq n-t$, and $\mu_s = 0$ if, and only if $s \in \{ 0,t,t'\}$. Then $P$ is of the form $P= z_0 v_0 + z_t v_t + z_{t'} v_{t'}$ for some $z_0,z_t,z_{t'} \in \C$.
Note that if $\gcd(t,n)= \gcd(t',n) = t_0 > 1$, then $P$ has at most $\frac{n}{t_0}$ pairwise distinct vertices; a contradiction. Thus, both $t$ and $t'$ are relatively prime to $n$. On the other hand, in this case the congruence $2t \equiv 0 \pmod{\frac{n}{\gcd(n,k)\gcd(n,m)}}$ implies that
$n \leq 2 \gcd(n,k) \gcd(n,m)$; a contradiction.

If the conditions in (b) hold, then the solutions are of the form $P= z_0 v_0 + z_t v_t$ for some $z_0, z_t \in \C$. Since the vertices of $P$
are pairwise distinct, it follows that $\gcd(t,n) = 1$, and thus, $P$ is a regular $n$-gon. Furthermore, as $m_1+m_2 \equiv k$ implies that $w \in \Re$, it follows that in this case $m_1+m_2 \not\equiv k \pmod n$.
If the conditions in (c) are satisfied, then the solutions are of the form $P= z_0 v_0 + z_t v_t + z_{n-t} v_{n-t}$ for some $z_0, z_t, z_{n-t} \in \C$. Similarly like in the previous case, it follows that $\gcd(t,n) = 1$. This condition implies that $P$ is an affinely regular $n$-gon. In addition, the congruence $m_1+m_2 \equiv k \pmod n$ follows from the fact that $w = \frac{\epsilon^{tm_1}- \epsilon^{tm_2}}{\epsilon^{tk}-1} \in \Re$.
This finishes the proof.

\section{Proof of Lemma~\ref{lem:conway}}\label{sec:conway}

We use the notation in the formulation of Lemma~\ref{lem:conway}. Observe that we may assume, without loss of generality, that $\gcd(k,k',l,l')$ and $n$ are relatively prime.

An elementary computation yields that for any value of $j$ with $1 \leq j \leq \frac{n}{2}$, we have $d_j = \sin \frac{j\pi}{n}$.
Consider values $k,l,k',l'$ such that $1 \neq \frac{\sin \frac{k\pi}{n}}{\sin \frac{l\pi}{n}} = \frac{\sin \frac{k'\pi}{n}}{\sin \frac{l'\pi}{n}}$.
By simple trigonometric identities, this equation can be transformed into
\begin{equation}\label{eq:cosines}
\cos \frac{(k-l')\pi}{n} - \cos \frac{(k+l')\pi}{n} - \cos \frac{(k'-l)\pi}{n} + \cos \frac{(k'+l)\pi}{n} = 0.
\end{equation}

To finish the proof we need the following theorem \cite[Theorem 7]{Conway}.

\begin{theorem}[Conway, Jones]\label{thm:conway}
Suppose we have at most four rational multiples of $\pi$ lying strictly between $0$ and $\frac{\pi}{2}$ for which some rational linear combination of their cosines is rational but no proper subset of them has this property. Then the appropriate linear combination is proportional to one from the following list.
\begin{itemize}
\item[(i)] $\cos \frac{\pi}{3} = \frac{1}{2}$.
\item[(ii)] $-\cos \phi + \cos \left( \frac{\pi}{3}-\phi\ \right) + \cos \left( \frac{\pi}{3} +\phi\right) = 0$, where $0 < \phi < \frac{\pi}{6}$.
\item[(iii)] $\cos \frac{\pi}{5} - \cos \frac{2\pi}{5} = \frac{1}{2}$.
\item[(iv)] $\cos \frac{\pi}{7}-\cos \frac{2\pi}{7}+\cos \frac{3\pi}{7} = \frac{1}{2}$.
\item[(vi)] $\cos \frac{\pi}{5}-\cos \frac{\pi}{15} + \cos \frac{4\pi}{15} = \frac{1}{2}$.
\item[(vii)] $- \cos \frac{2\pi}{5} + \cos \frac{2\pi}{15} - \cos \frac{7\pi}{15} = \frac{1}{2}$.
\item[(viii)] $\cos \frac{\pi}{7} + \cos \frac{3\pi}{7}-\cos \frac{\pi}{21} + \cos \frac{8\pi}{21} = \frac{1}{2}$.
\item[(ix)] $\cos \frac{\pi}{7}-\cos\frac{2\pi}{7}+\cos \frac{2\pi}{21}-\cos \frac{5\pi}{21} = \frac{1}{2}$.
\item[(x)] $-\cos \frac{2\pi}{7} + \cos \frac{3\pi}{7} + \cos \frac{4\pi}{21} + \cos \frac{10\pi}{21} = \frac{1}{2}$.
\item[(xi)] $-\cos \frac{\pi}{15} + \cos \frac{2\pi}{15} + \cos \frac{4\pi}{15} - \cos \frac{7\pi}{15} = \frac{1}{2}$.
\end{itemize}
 \end{theorem}

First, we consider the case that each of the four angles in (\ref{eq:cosines}) is an integer multiple of $\frac{\pi}{21}$.
This means that $\frac{21(k-l')}{n}, \frac{21(k+l')}{n},\frac{21(k'-l)}{n},\frac{21(k'+l)}{n}$ are all integers, implying that $n$ is a divisor of each of $42k,42k',42l,42l'$. Since $\gcd(k,k',l,l')$ and $n$ are relatively prime, from this it follows that $n$ is a divisor of $42$. On the other hand, an elementary computation shows that Lemma~\ref{lem:conway} holds for the lengths of the sides and diagonals of a regular $42$-gon. This shows the assertion in this case.
A similar argument proves the lemma in the case that each of the four angles in (\ref{eq:cosines}) is an integer multiple of $\frac{\pi}{15}$.

From now on, we assume that at least one of the angles in (\ref{eq:cosines}) is an integer multiple of neither $\frac{\pi}{21}$ nor $\frac{\pi}{15}$. First, we consider the case that no angle in (\ref{eq:cosines}) is an integer multiple of $\frac{\pi}{2}$, and without loss of generality, we assume that $ 0 < l' \leq l \leq k' \leq k \leq \frac{n}{2}$, which implies that $0 \leq k-l' \leq k'-l < \frac{n}{2}$, and $0 < k+l', k'+l \leq n$.
Using the trigonometric identities $\cos \alpha = \cos(-\alpha) = - \cos (\pi-\alpha)$ for $\alpha \in \Re$, it is easy to see that by Theorem~\ref{thm:conway}, two pairs of cosines cancel out in (\ref{eq:cosines}).
Note that since no member is zero, there are exactly two positive and two negative members, and that $\cos \frac{(k-l')\pi}{n}$ is positive and
$- \cos \frac{(k'-l)\pi}{n}$ is negative. Thus, we have three possibilities:
\begin{itemize}
\item $k-l'=k+l'$ and $k'-l=k'+l$. In this case $l=l'=0$; a contradiction.
\item $k-l'=k'-l$ and $k+l'=k'+l$. In this case $k=k'$ and $l=l'$, and the assertion follows.
\item $k-l'=n-k'-l$ and $k'-l=n-k-l'$. In this case $l=l'$ and $k+k'=n$, implying $k=k'=\frac{n}{2}$; a contradiction.
\end{itemize}
In the case that at least one of the angles is an integer multiple of $\frac{\pi}{2}$, Lemma~\ref{lem:conway} can be proven using a similar technique and case analysis.

\section{Symmetric realization of polytopes}\label{sec:polytopes}

In the proof for any $X \subset \Re^d$, we denote the affine hull of $X$ by $\aff X$, we let $\B^d$ be the closed unit ball with the origin as its center, and set $\Sph^{d-1} = \bd \B^d$.

First, note that (\ref{thm:symmetric}.2) or (\ref{thm:symmetric}.3) clearly implies (\ref{thm:symmetric}.1).
We prove that (\ref{thm:symmetric}.1) yields (\ref{thm:symmetric}.2).

As a first step, we show that the points $p_1,\ldots,p_{d+1}$ are affinely independent; in particular, we show, by induction on $s$, that for any $2 \leq s \leq d+1$, $\aff \{ p_1,p\ldots,p_s \}$ is an $(s-1)$-flat.
First, if $p_1 = p_2$, then $P$ is a single point, a contradiction, and thus, the statement holds for $s=2$. Now we assume that $\aff \{ p_1,\ldots,p_s \} $ is an $(s-1)$-flat for some $2 \leq s \leq d$, and show that $\aff \{ p_1,\ldots,p_{s+1} \}$ is an $s$-flat. Observe that by (\ref{thm:symmetric}.1), for every integer $j$, $\aff \{ p_{j+1},\ldots,p_{j+s} \}$ is also an $(s-1)$-flat. On the other hand, if $\aff \{ p_1,\ldots,p_{s+1} \}$ is not an $s$-flat, then $\aff \{ p_1,\ldots,p_{s+1} \} = \aff \{ p_1,\ldots,p_{s} \} = \aff \{ p_2,\ldots,p_{s+1} \}$, which, by (\ref{thm:symmetric}.1) yields that $P \subset \aff \{ p_1,\ldots,p_{s} \}$, a contradiction. Thus, we have that $p_1,\ldots,p_{d+1}$ are affinely independent.

Let $\phi : \Re^d \to \Re^d$ be the affine transformation defined by $\phi(p_s) = p_{s+1}$ for $s=1,2,\ldots,d+1$. Since $\conv \{ p_1,\ldots,p_{d+1}\}$ and $\conv \{ p_2,\ldots, p_{d+2}\}$ are congruent, $\phi$ is a congruence. Note that as $p_1,\ldots,p_{d+1}$ are affinely independent, for any $q \in \Re^d$, the distances of $q$ from these points determine $q$.
Thus, for any integer $j$, we have $\phi(p_j)=p_{j+1}$, and (\ref{thm:symmetric}.2) holds.

Finally, we prove that (\ref{thm:symmetric}.2) yields (\ref{thm:symmetric}.3).
Without loss of generality, let $\B^d$ be the unique smallest ball that contains $P$. Then $\B^d$ is the smallest ball containing $\phi(P)$ as well. Thus, $\phi$ is an isometry of $\B^d$, from which it follows that $p_j \in \Sph^{d-1}$ if, and only if $p_{j+1} \in \Sph^{d-1}$. This implies that $P$ is inscribed in $\Sph^{d-1}$.

We present two different arguments that finish the proof from this point.

\emph{First proof of (\ref{thm:symmetric}.3)}.
Let $E$ be the unique smallest volume ellipsoid containing $P$. Since $E$ is unique, $\Aut(P) \leq \Aut(E)$. On the other hand, the only ellipsoids whose symmetry groups contain an element $\phi$ of order $n \geq 5$ such that for some $p \in \Re^d$, the affine hull of the orbit of $p$ is $\Re^d$, are balls. Thus, without loss of generality, we may assume that $E=\B^d$.
We use the following, well-known properties of the smallest volume ellipsoid circumscribed about $P$ (cf. e.g. \cite{henk}).

\begin{theorem}\label{thm:Lowner}
Let $K\subset \B^d$ be a compact, convex set. Then $\B^d$ is the smallest volume ellipsoid circumscribed about $K$ if, and only if
for some $d\leq n \leq \frac{d(d+3)}{2}$ and $k = 1, \ldots , n$, there are $u_k\in \Sph^{d-1} \cap \bd K $ and $\lambda_k > 0$
such that
\begin{equation}\label{eq:lownerjohn}
0=\sum\limits_{k=1}^{n}\lambda_k u_k, \quad \Id =\sum\limits_{k=1}^n\lambda_k u_k\otimes u_k,
\end{equation}
where $\Id$ is the $d$-dimensional identity matrix, and for $u, v \in \Re^d$, $u\otimes v$ denotes the $d\times d$ matrix $u v^T$.
\end{theorem}

Thus, since $P \cap \B^d$ is the vertex set of $P$, there are some \emph{nonnegative} coefficients $\lambda_1,\lambda_2,\ldots,\lambda_n$
which, together with the points $p_1,p_2,\ldots,p_n$, satisfy the conditions in (\ref{eq:lownerjohn}).
We show that these points, with the coefficients $\lambda_1 = \lambda_2 = \ldots = \lambda_n = \frac{d}{n}$ also satisfy the conditions in (\ref{eq:lownerjohn}).
Indeed, set $0  < \lambda = \frac{\sum_{j=1}^n \lambda_j}{n}$. By (\ref{thm:symmetric}.2), we have that $\sum_{j=1}^n \lambda_{j+k} p_j = 0$ for every integer $k$, and thus, $\sum_{j+1}^n \lambda p_j = 0$. Similarly, since $\Id = \sum_{j=1}^n \lambda_{j+k} p_j \otimes p_j$ holds for for every integer $k$, it follows that $\Id = \lambda \sum_{j=1}^n p_j \otimes p_j$. Since $|p_j| = 1$ for every value of $j$, this equality implies that $d = \tr(\Id)= n \lambda$, that is, $\lambda = \frac{d}{n}$.

In the following we set $\bar{p}_k= \sqrt{\frac{d}{n}} p_k$ for $k=1,2,\ldots,n$, and observe that
\begin{equation}\label{eq:lownerjohn2}
0=\sum\limits_{k=1}^{n} \bar{p}_k, \quad \Id =\sum\limits_{k=1}^n  \bar{p}_k\otimes \bar{p}_k,
\end{equation}
Let $G$ be the Gram matrix of the vectors $\bar{p}_k$, that is, $G_{jk}= \langle \bar{p}_j, \bar{p}_k\rangle$. Then, using an elementary algebraic transformation, from the second equality in (\ref{eq:lownerjohn2}) we obtain that $G^2 = G$, and thus, that $G$ is the matrix of an orthogonal projection in $\Re^n$ into a $d$-dimensional subspace, with rank $d$.
This yields that $G$ has two eigenvalues, $1$ and $0$, with multiplicities $d$ and $n-d$, respectively, and, furthermore, we can write $G$ in the form $AA^T$, where $A$ is an $(n \times d)$ matrix, and the columns of $A$ form an orthonormal system in $\Re^d$. Equivalently, $G$ can be written in the form $G=B D B^T$, where $D$ is the diagonal matrix in which the first $d$ diagonal elements are equal to $1$, and the last $n-d$ elements are equal to $0$, and $B$ is an orthogonal matrix in $\Re^n$. Here, $A$ is the matrix composed of the first $d$ columns of $B$.
Observe that $v_j$ is the $j$th row of $A$, or equivalently, if we extend $A$ to an orthogonal matrix $B$, then the coordinates of $v_j$
are the first $d$ coordinates of the $j$th vector in the orthonormal system formed by the rows of $B$.

By (\ref{thm:symmetric}.2), $G$ is a circulant matrix.
Let $(c_0,c_1,\ldots,c_{n-1})$ be the first row of $G$. Then the eigenvalues of $G$ are $\mu_k = \sum_{j=0}^{n-1} c_j \epsilon^{jk}$, where $\epsilon = \cos \frac{2\pi}{n} + i \sin \frac{2\pi}{n}$, and the corresponding eigenvectors are $v_k = \frac{1}{\sqrt{n}} (1,\epsilon^k, \ldots, \epsilon^{(n-1)k} )$.
Now all eigenvalues are $0$ or $1$, and hence, $\mu_k = 1$ if and only if $\mu_{n-k}=1$.
Thus, the real $2$-flat spanned by $u_k= \frac{1}{2}(v_k+v_{n-k})$ and $u'_k = \frac{1}{2i}(v_k-v_{n-k})$ is contained in one of the two eigenspaces.
Let $F$ be the eigenspace associated to $1$. Then $F^{\perp}$ is the eigenspace associated to $0$.
The definition of $G$ and the fact that $\sum_{k=1}^n \bar{p}_k = 0$ yield that $\sum_{k=1}^n c_k = 0$. Thus, $\mu_0 = 0$.

If $n$ is even; that is, if $\mu_{n/2}$ exists, then $n-d$ is even if and only if $d$ is even. Hence, $\mu_0 = 0$ implies that if $n$ is even, then $\mu_{n/2} = 0$ if and only if $d$ is even. This yields that if $d$ is even, then $F$ is spanned by pairs of vectors of the form $u_k, u'_k$, where $0 < k < \frac{n}{2}$, and if $d$ is odd, then $F$ is spanned by the vector $(1,-1,\ldots, (-1)^n)$, and by pairs of vectors of the form $u_k, u'_k$, where $0 < k < \frac{n}{2}$. This shows that for some values $0 < k_1 < k_2 < \ldots < k_{\lfloor d/2 \rfloor} < \frac{n}{2}$, the Gram matrix of the vertices of $\sqrt{\frac{d}{n}} Q(k_1,k_2,\ldots,k_{\lfloor d/2 \rfloor})$ is equal to $G$. On the other hand, $G$ determines the pairwise distances between the points $\bar{p}_1, \bar{p}_2, \ldots, \bar{p}_n$. Thus, $P$ is congruent to $\sqrt{\frac{d}{n}} Q(k_1,k_2,\ldots,k_{\lfloor d/2 \rfloor})$, and the assertion readily follows.

\emph{Second proof of (\ref{thm:symmetric}.3)}.
Using the fact that $P$ is inscribed in $\Sph^{d-1}$ it follows that $\phi$ in (\ref{thm:symmetric}.2) is an orthogonal linear transformation. Let its matrix be denoted by $A$. Since $\phi$ is invertible, the diagonal elements in the Jordan form of $A$ are nonzero, from which it easily follows that $A$ is diagonalizable over $\C$. As $\phi$ is a real matrix, its complex eigenvalues are either real, or pairs of conjugate nonreal complex numbers. Note that if $z=x+iy \in \C^d$ is an eigenvector of $A$ associated to the eigenvalue $\lambda \notin \Re$, then the linear subspace in $\Re^d$, spanned by $x$ and $y$, is invariant under $\phi$. Thus, in a suitable orthonormal basis, the matrix of $\phi$ is a block diagonal matrix, where each block is either $1 \times 1$ (belonging to a real eigenvalue), or $2\times 2$ (belonging to a pair of conjugate complex eigenvalues). Let these blocks be $B_1, B_2, \ldots, B_k$.

As $\phi$ is orthogonal, each $1 \times 1$ block is either $1$ or $-1$ (corresponding to the identity, and to the reflection about the origin, respectively), and since $2\times 2$ blocks belong to nonreal eigenvalues (that is, they are not diagonalizable as $2\times 2$ matrices), they are $2$-dimensional rotation matrices. From the fact that the \emph{affine hull} of the vectors $p_1, \ldots, p_n$ is $\Re^d$, it follows that there is no block equal to $1$, and there at most one block equal to $-1$. Thus, if $d$ is even, then each block is a $2$-dimensional rotation, and if $d$ is odd, then one block is reflection about the origin, and every other block is a $2$-dimensional rotation. In the latter case we can assume that the last block belongs to the reflection, that is, $B_{(d+1)/2} = -1$.

Clearly, without loss of generality, we may assume that each angle of rotation is strictly less than $\pi$.
Thus, and since the order of each rotation is a divisor of $n$, for $i = 1,2, \ldots, \left\lfloor d/2 \right\rfloor$, we have
\[
B_i = 
\left[
\begin{array}{cc}
\cos \frac{2k_i \pi}{n} & -\sin \frac{2k_i \pi}{n} \\
\sin \frac{2k_i \pi}{n} & \cos \frac{2k_i \pi}{n}
\end{array}
\right]
\]
for some $0 < k_i < \frac{n}{2}$.
The fact that the affine hull of the points $p_1, p_2,\ldots, p_n$ is $\Re^d$ implies that the values $k_i$ are pairwise different.
Hence, setting $p_n = \frac{1}{\sqrt{\left\lfloor \frac{d+1}{2} \right\rfloor}} (1,0,1,0,\ldots )$, which we can do without loss of generality, and observing that $p_i = \phi^i (p_n)$, the assertion follows.

\section{Concluding remarks and questions}\label{sec:remarks}

\begin{remark}
Let $n,m_1,m_2,k \in \Z$ be integers such that $n \geq 4$,  neither $k$ nor $m_1-m_2$ is divisible by $n$, and $\gcd(n,k,m_1-m_2)=t > 1$. It is easy to see that there is an $n$-gon $P$, which is not affinely regular, but whose vertices satisfy the equations $p_{j+m_1}-p_{j+m_2} = w (p_{j+k}-p_j)$ for $j=1,2,\ldots,n$ and for some $w \in \C$ with $|w| \neq 1$ and independent of $j$. Indeed, let $P$ be a regular $n$-gon. Then $P$ satisfies the conditions for some suitable value of $w$. For $j=1,2,\ldots, t$, let $P_j$ denote the set of vertices of $P$ whose vertices are congruent to $j$ $\pmod \frac{n}{t}$. Then, translating $P_1, P_2, \ldots, P_t$ by any vectors $v_1,v_2,\ldots,v_t$, respectively, we obtain a polygon $P'$ which still satisfies our conditions.
\end{remark}

We note that the proof of Theorem~\ref{thm:main} describes exactly which are the values of $n,m_1,m_2,k$ such that if $P = (p_1,p_2,\ldots,p_n) \in \C^n$ satisfies the equations $p_{j+m_1}-p_{j+m_2} = w (p_{j+k}-p_j)$ for $j=1,2,\ldots,n$ and for some $w \in \C$ with $|w| \neq 1$, then $P$ is affinely regular.
In particular, the following holds.

\begin{remark}\label{rem:counterexamples}
Let $n,m_1,m_2,k \in \Z$ be integers such that $n \geq 4$,  neither $k$ nor $m_1-m_2$ is divisible by $n$, and $\gcd(n,k,m_1-m_2) = 1$.
Then there is an $n$-gon $P=(p_1,p_2,\ldots,p_n) \in \C^n$ which is not affinely regular and its vertices satisfy the equations $p_{j+m_1}-p_{j+m_2} = w (p_{j+k}-p_j)$ for $j=1,2,\ldots,n$ and for some $w \in \C$ with $|w| \neq 1$, if, and only if $n$ is even, and at least one of the following holds.
\begin{itemize}
\item[(i)] the congruence system $tk \equiv t'k \pmod n$, $tm_1+\frac{n}{2} \equiv t'm_2 \pmod n$ and $tm_2+\frac{n}{2} \equiv t'm_1 \pmod n$ has a solution $t,t'$ satisfying $\gcd(n,t)=1$, and $t' \not\equiv \pm t \pmod n$.
\item[(ii)] the congruence system $tk \equiv -t'k \pmod n$, $tm_1+tk + \frac{n}{2} \equiv t'm_1$ and $tm_2+tk + \frac{n}{2} \equiv t'm_2 \pmod n$ has a solution $t,t'$ satisfying $\gcd(n,t)=1$, and $t' \not\equiv \pm t \pmod n$.
\end{itemize}
\end{remark}

\begin{proof}
The necessity of the conditions in (i) or (ii) was shown in the proof of Theorem~\ref{thm:main}. On the other hand, if the conditions in, say, (i) are satisfied, then the vertices of any $n$-gon of the form $P(x) = x v_t + (1-x) v_{t'}$, where $x \in [0,1]$ satisfy the required equations. Here the $n$-gons $v_t = \left( 1, \epsilon^t, \epsilon^{2t}, \ldots, \epsilon^{(n-1)t}\right)$ and $v_{t'} = \left( 1, \epsilon^{t'}, \epsilon^{2t'}, \ldots, \epsilon^{(n-1)t'}\right)$ are not affinely related. Thus, it is easy to see that for some $x \in (0,1)$, $P(x)$ is not affinely regular.
\end{proof}

\begin{figure}[ht]
\includegraphics[width=0.8\textwidth]{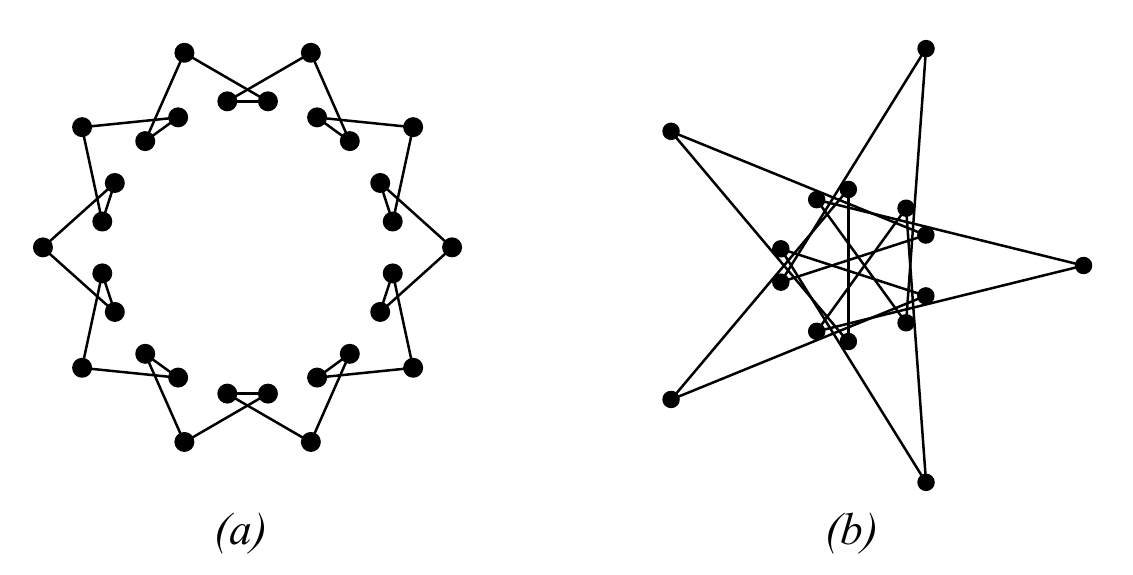}
\caption[]{Examples of polygons satisfying the conditions in Problem~\ref{prob:main}}
\label{fig:examples}
\end{figure}

Figure~\ref{fig:examples} shows examples for $n$-gons, which are not affinely regular, but satisfy the conditions in Problem~\ref{prob:main} for some values of $k,m_1,m_2 \in \Z$, and $w \in \C$. Panel (a) shows the $(30)$-gon $P= 0.8 v_1+0.2 v_{11}$. The vertices of $P$ satisfy the equations $p_{j+7}-p_{j+2} = w (p_{j+6}-p_j)$ for every index $j$, where $w = 0.809016\ldots + \mathrm{i} 0.2865\ldots$. Note that in this case $|w| = 0.850650\ldots \neq 1$, and $\gcd(30,6) \gcd(30,7-2) = 30$. Panel (b) shows the $(15)$-gon $P'=0.4 v_1+0.7 v_6 + 0.2 v_{11}$. The vertices of $P'$ satisfy the equations $p_{j+5} - p_{j+3} = w (p_{j+2}-p_j)$ for every index $j$, where $w = \cos \frac{2\pi}{5} + \mathrm{i} \sin \frac{2\pi}{5}$. Here $|w| = 1$.

\begin{remark}
We note that not all polytopes satisfying the conditions in Theorem~\ref{thm:symmetric} are cyclic. Indeed, let $d=4$, and for any $k \geq 3$, let $q_1,q_2,\ldots,q_k$ and $q'_1,q'_2,\ldots,q'_k$ be the vertices of two congruent regular $k$-gons, centered at the origin, where the first one is contained in the $(x_1,x_2)$ coordinate plane, and the second one in the $(x_3,x_4)$ coordinate plane. For $s=1,2,3,\ldots,k$, let $p_{2s}=q_s$, and $p_{2s+1}=q'_s$, and let $P=\conv \{ p_1,\ldots,p_{2k} \}$. Then $P$ satisfies the conditions in (\ref{thm:symmetric}.1), but it is not $2$-neighborly, and thus, it is not cyclic.
\end{remark}

\section{Acknowledgements}
The author expresses his gratitude to M. Nasz\'odi and K. Swanepoel for the fruitful discussions they had on this subject, to K. Swanepoel to direct his attention to the results in \cite{Conway}, and to an anonymous referee for his/her helpful comments, in particular to give an idea to prove (\ref{thm:symmetric}.3) in a different way.

\end{document}